\documentclass[11pt]{amsart}

\usepackage[T1]{fontenc}
\usepackage[latin1]{inputenc}
\usepackage{hyperref}
\usepackage{amscd}

\newcommand\Chow{{\rm Chow}}

\newtheorem{theorem}{Theorem}[section]

\newtheorem{corollary}[theorem]{Corollary} 
\newtheorem{proposition}[theorem]{Proposition} 
\newtheorem{question}[theorem]{Question} 
 
\newtheorem{definition}[theorem]{Definition}

\newtheorem{claim}[theorem]{Claim}

\newcommand\Pic{\mathop{\rm Pic}\nolimits}

\begin{document}

\title{A reduction map for nef line bundles} 

\date{\today}

\author[Bauer, Campana, Eckl, Kebekus et al]{Th.~Bauer, F.~Campana, 
  Th.~Eckl, S.~Kebekus, Th.~Peternell,  S.~Rams, T.~Szemberg, and L.~Wotzlaw}

\address{Thomas Bauer, Thomas Eckl, Stefan Kebekus and Thomas
  Peternell: Institut für Mathematik, Universität Bayreuth,
  95440~Bayreuth, Germany} 
\email{thomas.bauer@uni-bayreuth.de}
\email{thomas.eckl@uni-bayreuth.de}
\email{stefan.kebekus@uni-bayreuth.de}
\email{thomas.peternell@uni-bayreuth.de}

\address{Frédéric Campana: Département de Mathématiques, Université
  Nancy~1, BP~239, 54507~Vandoeuvre-les-Nance Cédex, France}
\email{campana@iecn.u-nancy.fr}

\address{Slawomir Rams: Mathematisches Institut der Universität,
  Bismarckstraße~$1\frac{1}{2}$, 91054~Erlangen, Germany }
\email{rams@mi.uni-erlangen.de}

\address{Tomasz Szemberg: Universität GH Essen, Fachbereich~6
  Mathematik, 45117~Essen, Germany} 
\email{mat905@uni-essen.de}

\address{Lorenz Wotzlaw: Mathematisches Institut, Humboldt-Universität Berlin, 
10099~Berlin, Germany}
\email{wotzlaw@mathematik.hu-berlin.de}

\maketitle
\tableofcontents

\section{Introduction}

In \cite{Ts00}, H.~Tsuji stated several very interesting assertions on
the structure of pseudo-effective line bundles $L$ on a projective
manifold $X$. In particular he postulated the existence of a
meromorphic ``reduction map'', which essentially says that through the
general point of $X$ there is a maximal irreducible $L$-flat
subvariety. Moreover the reduction map should be almost holomorphic,
i.e.~has compact fibers which do not meet the indeterminacy locus of
the reduction map. The proofs of \cite{Ts00} however are extremely
difficult to follow.

The purpose of this note is to establish the existence of a reduction
map in the case where $L$ is nef and to prove that it is almost
holomorphic ---this was also stated explicitly in \cite{Ts00}. Our
proof is completely algebraic while \cite{Ts00} works with deep
analytic methods. Finally, we show by a basic example that in the case
where $L$ is only pseudo-effective, the postulated reduction map
cannot be almost holomorphic ---in contrast to a claim in \cite{Ts00}.

These notes grew out from a small workshop held in Bayreuth in January
2001, supported by the Schwerpunktprogramm ``Global methods in complex
geometry'' of the Deutsche Forschungsgemeinschaft.

\section{A reduction map for nef line bundles}

In this section we want to prove the following structure theorem for
nef line bundles on a projective variety.

\begin{theorem}\label{thm:existence_of_reduction_map} 
  Let $L$ be a nef line bundle on a normal projective variety $X$.
  Then there exists an almost holomorphic, dominant rational map $f: X
  \dasharrow Y$ with connected fibers, called a ``reduction map'' such
  that
  \begin{enumerate} 
  \item $L$ is numerically trivial on all compact fibers $F$
    of $f$ with $\dim F = \dim X - \dim Y$
  \item for every general point $x \in X$ and every irreducible curve
    $C$ passing through $x$ with $\dim f(C) > 0$, we have $L \cdot C >
    0$.
  \end{enumerate}
  The map $f$ is unique up to birational equivalence of $Y$.
\end{theorem}

This theorem was stated without complete proof in Tsuji's paper
\cite{Ts00}. Relevant definitions are given now.

\begin{definition} 
  Let $X$ be an irreducible reduced projective complex space
  (projective variety, for short). A line bundle $L$ on $X$ is
  numerically trivial, if $L \cdot C = 0$ for all irreducible curves
  $C \subset X$. The line bundle $L$ is nef if $L \cdot C \geq 0$ for
  all curves $C$.
\end{definition}

Let $f: Y \to X$ be a surjective map from a projective variety $Y$.
Then clearly $L$ is numerically trivial (nef) if and only $f^*(L)$ is.

\begin{definition} 
  Let $X$ and $Y$ be normal projective varieties and $f: X \dasharrow
  Y$ a rational map and let $X^0 \subset X$ be the maximal open subset
  where $f$ is holomorphic.  The map $f$ is said to be almost
  holomorphic if some fibers of the restriction $f|_{X^0}$ are
  compact.
\end{definition}

\subsection{Construction of the reduction map}

\subsubsection{A criterion for numerical triviality}
\label{sec:criterion_ntriv}

In order to prove theorem~\ref{thm:existence_of_reduction_map} and
construct the reduction map, we will employ the following criterion
for a line bundle to be numerically trivial.

\begin{theorem} \label{thm:charact_of_nef}
  Let $X$ be an irreducible projective variety which is not
  necessarily normal.  Let $L$ be a nef line bundle on $X$. Then $L$
  is numerically trivial if and only if any two points in $X$ can be
  joined by a (connected) chain $C$ of curves such that $L \cdot C =
  0$.
\end{theorem}

In the remaining part of the present section~\ref{sec:criterion_ntriv}
we will prove theorem~\ref{thm:charact_of_nef}. The proof will be
performed by a reduction to the surface case. The argumentation is
then based on the following statement, which, in the smooth case, is a
simple corollary to the Hodge Index Theorem.

\begin{proposition}\label{prop:weak_criterion}
  Let $S$ be an irreducible, projective surface which is not
  necessarily normal, and let $q:S \to T$ be a morphism with connected
  fibers onto a curve. Assume that there exists a nef line bundle $L
  \in \Pic(S)$ and a curve $C \subset S$ such that $q(C)=T$ and such
  that
  $$
  L \cdot F = L \cdot C = 0
  $$
  holds, where $F$ is a general $q$-fiber. Then $L$ is numerically
  trivial.
\end{proposition}

\begin{proof}
  If $S$ is smooth, set $D = C + nF$, where $n$ is a large positive
  integer. Then we have $D^2 > 0$.  By the Hodge Index Theorem it
  follows that
  $$ 
  (L \cdot D)^2 \geq L^2 \cdot D^2 ,
  $$
  hence $L^2 = 0$, since by our assumptions $L \cdot D = 0$. So
  equality holds in the Index Theorem and therefore $L$ and $D$ are
  proportional: $L \equiv kD$ for some rational number $k$. Since $0 =
  L^2 = k^2 D^2$ and $D^2 > 0$, we conclude that $k = 0$. That ends
  the proof in the smooth case.
  
  If $S$ is singular, let $\delta: \tilde S \to S$ be a
  desingularisation of $S$ and let $\tilde C \subset \tilde S$ be a
  component of $\delta^{-1}(C)$ which maps surjectively onto $C$. Note
  that the fiber of $q\circ \delta$ need no longer be connected and
  consider the Stein factorisation
  $$
  \begin{CD}
    {\tilde S}      @>{\delta}>{\text{\scriptsize desing.}}> S \\
    @V{\tilde q}VV            @VV{q}V \\
    {\tilde T} @>>{\text{\scriptsize finite}}> T
  \end{CD}
  $$
  It follows immediately from the construction that $\tilde q
  (\tilde C)=\tilde T$, that $\delta^*(L)$ has degree 0 on $\tilde C$
  and on the general fiber of $\tilde q$. The argumentation above
  therefore yields that $\delta^*(L)$ is trivial on $\tilde S$. The
  claim follows.
\end{proof}

\subsubsection{Proof of Theorem~\ref{thm:charact_of_nef}}
If $L$ is numerically trivial, the assertion of
theorem~\ref{thm:charact_of_nef} is clear. We will therefore assume
that any two points can be connected by a curve which intersects $L$
with multiplicity 0, and we will show that $L$ is numerically trivial.
To this end, choose an arbitrary, irreducible curve $B \subset X$. We
are finished if we show that $L\cdot B=0$.

Let $a \in X$ be an arbitrary point which is not contained in $B$. For
any $b \in B$ we can find by assumption a connected, not necessarily
irreducible, curve $Z_b$ containing $a$ and $b$ such that $L \cdot Z_b
= 0$. Since the Chow variety has compact components and only a
countable number of components, we find a family $(Z_t)_{t\in T}$ of
curves, parametrised by a compact irreducible curve $T \subset
\Chow(X)$ such that for every point $b \in B$, there exists a point
$t\in T$ such that the curve $Z_t$ contains both $a$ and $b$. We
consider the universal family $S \subset X\times T$ over $T$ together
with the projection morphisms
$$
\begin{CD}
  S @>{p}>> X \\
  @V{q}VV \\
  T
\end{CD}
$$

\subsubsection*{Claim 1} 
There exists an irreducible component $S_0 \subset S$ such that
$p^*(L)$ is numerically trivial on $S_0$.

\subsubsection*{Proof of claim 1}
As all curves $Z_t$ contain the point $a$, the surface $S$ contains
the curve $\{a\} \times T$. Let $S_0 \subset S$ be a component which
contains $\{a\} \times T$. Since $\{a\} \times T$ intersects all
fibers of the natural projection morphism $q$, and since $p^*(L)$ is
trivial on $\{a\} \times T$, an application of
proposition~\ref{prop:weak_criterion} yields the claim. \qed

\subsubsection*{Claim 2} 
The bundle $p^*(L)$ is numerically trivial on $S$. 

\subsubsection*{Proof of claim 2}
We argue by contradiction and assume that there are components $S_j
\subset S$ where $p^*(L)$ is not numerically trivial. We can therefore
subdivide $S$ into two subvarieties as follows:
$$
\begin{array}{lcr}
S_{triv} & := & \{ \text{union of the irreducible components $S_i \subset 
S$}\quad \\
         &    & \quad \text{where $p^*(L)|_{S_i}$ is numerically trivial}\} \\
S_{ntriv} & := & \{ \text{union of the irreducible components $S_i \subset S$} 
\quad\\
         &    & \quad \text{where $p^*(L)|_{S_i}$ is not numerically trivial}\} 
\\
\end{array}
$$
By assumption, and by claim~1, both varieties are not empty. Since
$S$ is the universal family over a curve in $\Chow(X)$, the morphism
$q$ is equidimensional. In particular, since all components $S_i
\subset S$ are two-dimensional, every irreducible component $S_i$ maps
surjectively onto $T$. Thus, if $t \in T$ is a general point, the
connected fiber $q^{-1}(t)$ intersects both $S_{triv}$ and
$S_{ntriv}$. Thus, over a general point $t\in T$, there exists a point
in $S_{triv} \cap S_{ntriv}$. It is therefore possible to find a curve
$D \subset S_{triv} \cap S_{ntriv}$ which dominates $T$.

That, however, contradicts proposition~\ref{prop:weak_criterion}: on
one hand, since $D \subset S_{triv}$, the degree of $p^*(L)|_D$ is 0.
On the other hand, we can find an irreducible component $S_j \subset
S_{ntriv} \subset S$ which contains $D$. But because $p^*(L)$ has
degree 0 on the fibers of $p^*(L)|_{S_j}$,
proposition~\ref{prop:weak_criterion} asserts that $p^*(L)$ is
numerically trivial on $S_j$, contrary to our assumption.  This ends
the proof of Claim~2. \qed

We apply Claim~2 as follows: if $B' \subset S$ is any component of the
preimage $p^{-1}(B)$, then $p^*(L).B' = 0$. That shows that $L.B =0$,
and the proof of theorem~\ref{thm:charact_of_nef} is done. \qed

\subsubsection{Proof of Theorem~\ref{thm:existence_of_reduction_map}}
In order to derive Theorem~\ref{thm:existence_of_reduction_map} from
Theorem~\ref{thm:charact_of_nef}, we introduce an equivalence relation
on $X$ with setting $x \sim y$ if $x$ and $y$ can be joined by a
connected curve $C$ such that $L \cdot C = 0$. Then by \cite{Ca81} or
\cite[appendix]{Ca94}, there exists an almost holomorphic map $f: X
\dasharrow Y$ with connected fibers to a normal projective variety $Y$
such that, two general points $x$ and $y$ satisfy $x \sim y$ if and
only if $f(x) = f(y)$. This map $f$ gives the fibration we are looking
for.  

If $F$ is a general fiber, then $L \vert F \equiv 0$ by Theorem 2.4.

We still need to verify that $L \cdot C = 0$ for all curves $C$
contained in an {\it arbitrary} compact fiber $F_0$ of dimension $\dim
F_0 = \dim X - \dim Y$. To do that, let $H$ be an ample line bundle on
$X$ and pick
$$
D_1, \ldots, D_k \in \vert mH \vert 
$$
for $m $ large such that
$$ 
D_1 \cdot \ldots \cdot D_k \cdot F_0 = C+C' 
$$
with an effective curve $C'.$ Then
$$
L \cdot (C+C') = L \cdot D_1 \cdot \ldots \cdot D_k \cdot F
$$
with a general fiber $F$ of $f$, hence $L \cdot (C+C') = 0$. Since
$L$ is nef, we conclude $L \cdot C = 0.$ \qed

\subsection{Nef cohomolgy classes}

In Theorems~\ref{thm:existence_of_reduction_map} and
\ref{thm:charact_of_nef} we never really used the fact that $L$ is a
line bundle; only the property that $c_1(L) $ is a nef class is
important and even rationality of the class does not play any role.
Hence our results directly generalize to nef cohomology classes of
type $(1,1).$ To be precise, we fix a projective manifold (we stick to
the smooth case for sakes of simplicity) and we say that a class
$\alpha \in H^{1,1}(X,\mathbb R)$ is nef, if it is in the closure of
the cone generated by the K\"ahler classes. Moreover $\alpha$ is
numerically trivial, if $\alpha \cdot C = 0$ for all curves $C \subset
X$.

If $Z \subset X$ is a possibly singular subspace, then we say that
$\alpha $ is numerically trivial on $Z$, if for some (and hence for
all, see [Pa98]) desingularisation $\hat Z \to Z,$ the induced form
$f^*(\alpha)$ is numerically trivial on $\hat Z,$ i.e.  $f^*(\alpha)
\cdot C = 0$ for all curves $C \subset \hat Z.$ Here $f: \hat Z \to X$
denotes the canonical map. Similarly we define $\alpha$ to be nef on
$Z.$ If $Z$ is smooth, this is the same as to say that $\alpha \vert
Z$ is a nef cohomology class in the sense that $\alpha \vert Z$ is in
the closure of the K\"ahler cone of $Z$.

\begin{theorem} 
  Let $\alpha$ be a nef cohomology class on a smooth projective
  variety $X.$ Then there exists an almost holomorphic, dominant
  rational map $f: X \dasharrow Y$ with connected fibers, such that
  \begin{enumerate} 
  \item $\alpha$ is numerically trivial on all compact fibers $F$ of
    $f$ with $\dim F = \dim X - \dim Y$
  \item for every general point $x \in X$ and every irreducible curve
    $C$ passing through $x$ with $\dim f(C) > 0$, we have $\alpha
    \cdot C > 0$.
  \end{enumerate}
  The map $f$ is unique up to birational equivalence of $Y$.
  
  In particular, if two general points of $X$ can be joined by a chain
  $C$ of curves such that $\alpha \cdot C = 0,$ then $\alpha \equiv
  0$.
\end{theorem}

\subsection{The nef dimension}
\label{sec:construction}

Since $Y$ is unique up to a birational map, its dimension $\dim Y$ is
an invariant of $L$ which we compare to the other known invariants.

\begin{definition} 
  The dimension $\dim Y$ is called the nef dimension of $L$. We write
  $n(L) := \dim Y$.
\end{definition}

As usual we let $\nu(L)$ be the numerical Kodaira dimension of $L$,
i.e.~the maximal number $k$ such that $L^k \not \equiv 0$.
Alternatively, if $H$ is a fixed ample line bundle, then $\nu(L)$ is
the maximal number $k$ such that $L^k \cdot H^{n-k} \not = 0$.

\begin{proposition} \label{prop:numerical_kodaira_dim}
  The nef dimension is never smaller than the numerical Kodaira
  dimension:
  $$
  \nu (L) \leq n(L).
  $$
\end{proposition}
\begin{proof}
  Fix a very ample line bundle $H\in \Pic(X)$ and set $\nu := \nu
  (L)$. Let $Z$ be a general member cut out by $n-\nu$ elements of
  $\vert H \vert$. The dimension of $Z$ will thus be $\dim Z = \nu$,
  and since $L^{\nu} \cdot H^{n-\nu} > 0$, the restriction $L|_Z$ is
  big (and nef). Consequence: $\dim f(Z) = \nu$, since otherwise $Z$
  would be covered by curves $C$ which are contained in general fibers
  of $f$, so that $L \cdot C = 0$, contradicting the bigness of
  $L|_Z$. In particular, we have $\dim Y \geq \dim f(Z) = \nu$ and our
  claim is shown.
\end{proof}

\begin{corollary} 
  The nef dimension is never smaller than the Kodaira dimension:
  $$
  \kappa (L) \leq n(L).
  $$
  \qed
\end{corollary}

\subsection{The structure of the reduction map}

Again let $L$ be a nef line bundle on a projective manifold and $f : X
\to Y$ a reduction map. In this section we will shed light on the
structure of the map in a few simple cases. 

At the present time, however, we cannot say much about $f$ in good
generality. The following natural question is, to our knowledge, open.

\begin{question} 
  Are there any reasonable circumstances under which the reduction map
  can be taken holomorphic? Is it possible to construct an example
  where the reduction map cannot be chosen to be holomorphic?
\end{question}

Of course, the abundance conjecture implies that the $f$ can be chosen
holomorphic in the case where $L = K_X$. Likewise, we expect a
holomorphic reduction map when $L = -K_X$.

\subsubsection{The case where $L$ is big}
If $L$ is big and nef, then $n(L) = \dim X$. The converse however is
false: there are examples of surfaces $X$ carrying a line bundle $L$
such that $L \cdot C > 0$ for all curves $C$ but $L^2 = 0$. So in that
case $n(L) = 2,$ but $\nu (L) = 1$. The first explicit example is due
to Mumford, see \cite{Ha70}. This also shows that equality can fail in
proposition~\ref{prop:numerical_kodaira_dim}.

\subsubsection{The case where $n(L)=\dim X = 2$}
In this case, we have $L \cdot C > 0$ for all but an at most countable
number of irreducible curves $C \subset X$. These curves necessarily have 
semi-negative self-intersection $C^2
\leq 0$.

\subsubsection{The case where $n(X)=1$ and $\dim X = 2$ } \label{sec:n_is_1}

Here we can choose $Y$ to be a smooth curve. The almost-holomorphic
reduction $f:X\to Y$ will thus be holomorphic. In this situation, we
claim that the numerical class of a suitable multiple of $L$ comes
from downstairs:

\begin{proposition}
  If $n(X)=1$, then there exists a $\mathbb Q$-divisor $A$ on $Y$ such
  that
  $$ 
  L \equiv f^*(A).
  $$ 
\end{proposition}
\begin{proof}
  Since $L$ is nef, it carries a singular metric with positive
  curvature current $T,$ in particular $c_1(L) = [T].$ Since $L \cdot
  F = 0,$ and since every $(1,1)-$form $\varphi$ on $Y$ is closed, we
  conclude that
  $$
  f_*(T)(\varphi) = T(f^*(\varphi)) = [T] \cdot [f^*(\varphi)] =
  0,
  $$
  since $[f^*(\varphi)] = \lambda [F]$ for a fiber $F$ of $f$.
 
  Hence $f_*(T) = 0$.  Let $Y_0$ be the maximal open subset of $Y$
  over which $f$ is a submersion and let $X_0 = f^{-1}(Y_0)$. Then by
  \cite[(18)]{HL83}
  $$
  T \vert X_0 = \sum \mu_i F_i,
  $$
  where $F_i$ are fibers of $f \vert_{X_0}$ and $\mu_i > 0$. Hence
  ${\rm supp}(T - \sum \mu_i F_i) \subset X \setminus X_0$, an
  analytic set of dimension 1. Hence by a classical theorem (see
  e.g.~\cite{HL83}),
  $$
  T - \sum \mu_i F_i = \sum \lambda_j G_j
  $$
  where $G_j$ are irreducible components of $X \setminus X_0$,
  i.e.~of fibers, of dimension 1 and where $\lambda_j > 0.$ Since
  $(\sum \lambda_j G_j)^2 = 0$, Zariski's lemma shows that $\sum
  \lambda_j G_j$ is a multiple of fibers, whence $L \equiv f^*(A)$ for
  a suitable $\mathbb Q$-divisor $A$ on $B$.

Here is an algebraic proof, which however does not easily extend to
higher dimensions as the previous (see 2.3.4). \\
Fix an ample line bundle $A$ on $X$ and choose a positive rational
number $m$ such that $L \cdot A = mF \cdot A.$ Let $D = mF.$
Then
$L^2 = L \cdot D = D^2$ and $L \cdot A = D \cdot A.$ 
Introducing $r = L^2$ and $d = L \cdot A$, the ample line bundles $L' = L + A$
and $D' = D + A$ fulfill the equation
$$ (L' \cdot D')^2 = (r+2d+A^2)^2 = (L')^2 (D')^2,$$
hence $L$ and $D$ are numerically proportional, and so do $L$ and $D$. 
Consequently $L \equiv  f(A)$ as before. 
\end{proof}

\subsubsection{The case when $n(L)=1$ and $\dim X$ is arbitrary} 
The considerations of section~\ref{sec:n_is_1} easily generalize to
higher dimensions: if $X$ is a projective manifold and $L$ a nef line
bundle on $X$ with $n(L)=1,$ then the reduction map $f: X \to Y$ is
holomorphic and there exists a $\mathbb Q$-divisor on $Y$ such that $L
\equiv f^*(A).$

\subsubsection{The case where $n(L)=0$} 
We have $n(L) = 0$ if and only if $L \equiv 0$.

\section{A Counterexample}

In \cite{Ts00}, H.~Tsuji claims the following:
\begin{claim} \label{claim:tsuji}
  Let $(L,h)$ be a line bundle with a singular hermitian metric on a
  smooth projective variety $X$ such that the curvature current
  $\Theta_h$ is non-negative. Then there exists a (up to birational
  equivalence) unique rational fibration
  $$
  f: X --\rightarrow Y
  $$
  such that
  \begin{itemize}
  \item[(a)] $f$ is regular over the generic point of $Y$,
  \item[(b)] $(L,h)_{|F}$ is well defined and numerically trivial for
    every very general fiber~$F$, and
  \item[(c)] $\dim Y$ is minimal among such fibrations.
  \end{itemize} 
\end{claim}
Tsuji calls a pair $(L,h)$ numerically trivial if for every
irreducible curve $C$ on $X$, which is not contained in the singular
locus of $h$,
$$
(L,h).C = 0
$$
holds. The intersection number has the following
\begin{definition}
  Let $(L,h)$ be a line bundle with a singular hermitian metric on a
  smooth projective variety $X$ such that the curvature current
  $\Theta_h$ is non-negative. Let $C$ be an irreducible curve on $X$
  such that the natural morphism $\mathcal{I}(h^m) \otimes
  \mathcal{O}_C \rightarrow \mathcal{O}_C$ is an isomorphism at the
  generic point of $C$ for every $m \geq 0$.  Then
  $$
  (L,h).C := \limsup_{m \rightarrow \infty} 
  \frac{\dim H^0(C, \mathcal{O}_C(mL) \otimes \mathcal{I}(h^m)/
  \mathrm{tor})}{m} ,
  $$
  where $\mathrm{tor}$ denotes the torsion part of
  $\mathcal{O}_C(mL) \otimes \mathcal{I}(h^m)$.
\end{definition}
Of course, claim~\ref{claim:tsuji} is trivial as stated ---except for
the uniqueness assertion. What is meant is that every curve $C$ with
$(L,h).C = 0$ is contracted by $f$. This is also clear from Tsuji's
construction.

In order to make sense of claim~\ref{claim:tsuji}, statement (c) must
be given the following meaning: if $x \in X$ is very general, then
$(L,h) \cdot C > 0$ for all curves $C$ through $x$.

There is an easy counterexample to the Claim~\ref{claim:tsuji}: Take
$X = \mathbb P^2$ with homogeneous coordinates $(z_0:z_1:z_2)$, $L =
\mathcal{O}(1)$ and let $h$ be induced by the incomplete linear system
of lines passing through $(1:0:0)$. Then the corresponding
plurisubharmonic function $\phi$ is given by $\frac{1}{2}\log(|z_1|^2
+ |z_2|^2)$.

By \cite[5.9]{De00}, it is possible to reduce the calculation of
$\mathcal{I}(h^m)$ to an algebraic problem: Since the ideal sheaf
$\mathcal{J}_p$ generated by $z_1, z_2$ describes the reduced point $p
= (1:0:0)$, it is integrally closed. Let $\mu: \mathbb{F}_1
\rightarrow \mathbb{P}^2$ be the blow up of $\mathbb{P}^2$ in $p$.
Then $\mu^\ast \mathcal{J}_p$ is the invertible sheaf
$\mathcal{O}(-E)$ associated with the exceptional divisor $E$.  Now,
one has $K_{\mathbb{F}_1} = \mu^\ast K_{\mathbb{P}^2} + E$.  By the
direct image formula in \cite[5.8]{De00} it follows
$$
\mathcal{I}(m \phi) = \mu_\ast(
\mathcal{O}(K_{\mathbb{F}_1} - \mu^\ast K_{\mathbb{P}^2}) \otimes
\mathcal{I}(m \phi \circ \mu)).
$$
Now, $(z_i \circ \mu)$ are generators of the ideal
$\mathcal{O}(-E)$, hence
$$
m \phi \circ \mu \sim m \log g,
$$
where $g$ is a local generator of $\mathcal{O}(-E)$.  But
$\mathcal{I}(m \phi \circ \mu)) = \mathcal{O}(-mE)$, hence
$$
\mathcal{I}(m \phi) = \mu_\ast \mathcal{O}_{\mathbb{F}_1}((1-m)E) =
\mathcal{J}_p^{m-1}.
$$
   
Let $C \subset \mathbb{P}^2$ be a smooth curve of degree $d$ and genus
$g$.  If $p \not\in C$ then $\mathcal{O}_C(mL) \otimes
\mathcal{J}_p^{m-1} = \mathcal{O}_C(mL)$, and $(L,h).C = L.C = d$.  On
the other hand, if $p \in C$, then
$$
\mathcal{O}_C(mL) \otimes \mathcal{J}_p^{m-1} = 
\mathcal{O}_C(mL - (m-1)p).
$$
The degree of this invertible sheaf is $md - m + 1 = m(d-1) + 1$,
hence by Riemann-Roch
\begin{eqnarray*}
  \limsup_{m \rightarrow \infty} 
  \frac{\dim H^0(C, \mathcal{O}_C(mL) \otimes \mathcal{I}(m\phi))}{m} 
  & = &
  \limsup_{m \rightarrow \infty} \frac{m(d-1) + 1 + 1 -g}{m} \\ 
  & = & d-1. 
\end{eqnarray*}
   
Consequently, Tsuji's fibration could only be the trivial map
$\mathbb{P}^2 \rightarrow \mathbb{P}^2$: It can't be the map to a
point, because there are curves with intersection number $\geq 1$ on
$\mathbb{P}^2$, and it can't be an almost holomorphic rational map to
a curve, because $\mathbb P_2$ does not contain curves with vanishing
self-intersection.  This contradicts condition (c) in the claim.
   
It remains an open question if claim~\ref{claim:tsuji} is true without
condition (a), and if in the case of singular hermitian metrics
induced by linear systems, the fibration may be taken as the induced
rational map.


\begin{thebibliography}{HL83}

  
\bibitem[Ca81]{Ca81} F.~Campana. Cor\'eduction alg\'ebrique d'un
  espace analytique faiblement k\"ahl\'erien compact. {\em Inv.
    Math.}, 63:187-223, 1981.

\bibitem[Ca94]{Ca94} F.~Campana. Remarques sur le revetement universel des 
vari\'eti\'es
k\"ahleriennes compactes. {\em Bull. SMF}, 122:255-284, 1994

\bibitem[De00]{De00} J.P.~Demailly. Vanishing theorems and effective
  results in algebraic geometry. Trieste lecture notes, 2000

\bibitem[Ha70]{Ha70} R.~Hartshorne. {\em Ample subvarieties of
    algebraic varieties}, volume 156 of {\em Lecture Notes in
    Mathematics}, Springer, 1970.
  
\bibitem[HL83]{HL83} R.~Harvey, H.B.~Lawson. An intrinsic
  characterization of K\"ahler manifolds. {\em Inv. math.},
  74:169--198, 1983.

\bibitem[Pa98]{Pa98} M.Paun. Sur l'effectivit\'e num\'erique des images inverses
de fibr\'es en droites. {\em Math. Ann. }, 310:411--421, 1998.
  
\bibitem[Ts00]{Ts00} H.~Tsuji. Numerically trivial fibrations.
  LANL-preprint math.AG/0001023, October 2000.

\end{thebibliography}
\end{document}